\documentclass[11pt]{article}
\usepackage{amscd}
\usepackage{amsfonts}
\usepackage{amsmath}
\usepackage{amssymb}
\usepackage{amsthm}
\usepackage{bbm}
\usepackage{CJK}
\usepackage{fancyhdr}
\usepackage{graphicx}
\usepackage{hyperref}
\usepackage{indentfirst}
\usepackage{latexsym}
\usepackage{mathrsfs}
\usepackage{xypic}

\usepackage[top=1in,bottom=1in,left=1.25in,right=1.25in]{geometry}
\newtheorem{theorem}{Theorem}[section]

\newtheorem{definition}[theorem]{Definition}
\newtheorem{proposition}[theorem]{Proposition}
\newtheorem{example}[theorem]{Example}

\usepackage[top=1in,bottom=1in,left=1.25in,right=1.25in]{geometry}
\def\<{\langle}
\def\>{\rangle}
\def\a{\alpha}
\def\b{\beta}

\def\c{\cdot}

\def\g{\gamma}

\def\lr{\longrightarrow}

\def\o{\otimes}

\date{}
\begin{document}
\renewcommand{\baselinestretch}{1.2}
\renewcommand{\arraystretch}{1.0}
\title{\bf Manin triples and quasitriangular structures  of Hom-Poisson bialgebras}
\author{{\bf Shuangjian Guo$^{1}$, Xiaohui Zhang$^{2}$,  Shengxiang Wang$^{3}$\footnote
        { Corresponding author(Shengxiang Wang):~~wangsx-math@163.com} }\\
{\small 1. School of Mathematics and Statistics, Guizhou University of Finance and Economics} \\
{\small  Guiyang  550025, P. R. of China} \\
{\small 2.  School of Mathematical Sciences, Qufu Normal University}\\
{\small Qufu  273165, P. R. of China}\\
{\small 3.~ School of Mathematics and Finance, Chuzhou University}\\
 {\small   Chuzhou 239000,  P. R. of China}}
 \maketitle

\begin{center}
\begin{minipage}{13.cm}

{\bf \begin{center} ABSTRACT \end{center}}
In this paper, we first introduce the definition of a Hom-Poisson bialgebra
and give an equivalent descriptions via the Manin triple of Hom-Poisson algebras.
Also we introduce notions of $\mathcal{O}$-operator on a Hom-Poisson algebra,
post-Hom-Poisson algebra  and quasitriangular Hom-Poisson bialgebra, and present a method to
 construct post-Hom-Poisson algebras.
Finally, we show that a quasitriangular Hom-Poisson bialgebra naturally yields a post-Hom-Poisson algebra.

{\bf Key words}: Hom-Poisson algebras, Hom-Poisson bialgebras, Manin triple, $r$-matrics, $\mathcal{O}$-operators.

{\bf 2010 Mathematics Subject Classification:} 16T05; 17A30; 17B62.
 \end{minipage}
 \end{center}
 \normalsize\vskip1cm

\section*{Introduction}
A Poisson algebra is both a Lie algebra and a commutative associative algebra which are
compatible in a certain sense. Poisson algebras play important roles in many fields in mathematics
and mathematical physics, such as the Poisson geometry, integrable systems, non-commutative
(algebraic or differential) geometry, and so on (see \cite{V94} and the references therein).

As generalizations of Lie algebras, Hom-Lie algebras were introduced motivated by
applications to physics and to deformations of Lie algebras, especially Lie algebras
of vector fields. The notion of Hom-Lie algebras was firstly introduced by Hartwig,
Larsson and Silvestrov  to describe the structure of certain $q$-deformations
of the Witt and the Virasoro algebras, see \cite{HLS06}.
More precisely, a  Hom-Lie algebras are different from
Lie algebras as the Jacobi identity is replaced by a twisted form using a morphism.
This twisted Jacobi identity is called Hom-Jacobi identity given by
\begin{eqnarray*}
[\alpha(x),[y,z]]+[\alpha(y),[z,x]]+[\alpha(z),[x,y]]=0.
\end{eqnarray*}

The twisting parts of the defining identities were transferred to other algebraic structures.
In \cite{MS08},  Makhlouf and  Silvestrov introduced the notion of
 Hom-associative algebras.
  In  \cite{MY14}, Makhlouf and Yau introduced the definition of Hom-Poisson algebras,  which is a natural
generalization of  Poisson  algebras, Hom-Lie algebra and Hom associative algebras.
Since then it was further investigated by \cite{BEM12}. Sheng and Bai introduced a new definition of a Hom-Lie bialgebra and investigated their
properties by \cite{SB14}.

In this paper, we introduce a definition  of a Hom-Poisson algebra. Now given a Hom-Poisson bialgebra $P$, $P\oplus P^{\ast}$ is a Hom-Poisson algebra such that $(P\oplus P^{\ast}, P, P^{\ast})$ is a Manin triple of Hom-Poisson. We also study that a Hom-Poisson bialgebra exhibits both the features of a Hom-Lie bialgebra and an infinitesimal
Hom-bialgebra. In particular, a Hom-Poisson bialgebra can be obtained by solving the classical Hom-Yang-Baxter
equation (CHYBE) and Hom-associative Yang-Baxter equation (HAYBE) uniformly (namely, Hom-Poisson Yang-
Baxter equation (HPYBE)). In detail, and construct an $r$-matrix in the semirdirct Hom-Poisson algebra by introducing a notation of an $\mathcal{O}$ operator for a Hom-Poisson algebra.

The paper is organized as follows. In Section 2, we introduce a notion of matched pair of Hom-Poisson
algebras which plays an essential role. In Section 3, we introduce notions of Manin triple
of Hom-Poisson algebras and Hom-Poisson bialgebras and then give the equivalence between them in terms of
the matched pairs of Hom-Poisson algebras. In Section 4, we consider the coboundary cases and get PYBE
which is a combination of CYBE and AYBE. In Section 5, we introduce notions of $\mathcal{O}$-operator of
weight $\lambda \in k$ of a Hom-Poisson algebra and post-Hom-Poisson algebra and interpret the close relationships
between them and Hom-Poisson bialgebras.

\section{Preliminary}

\def\theequation{1.\arabic{equation}}
\setcounter{equation} {0} \hskip\parindent

In this paper, the base field is taken to be $k$ unless otherwise specified. This is the field over
which we take all the associative, Lie and Poisson algebras, vector spaces, linear maps and tensor
products, etc. and the characteristic of $k$ is zero.

In this section, we will recall from  \cite{BM14} and \cite{MS08}   that the basic definitions and results
 on the  Hom-associative algebras and  Hom-Lie algebras.

\begin{definition}

 A Hom-associative algebra is a triple $(A,\mu,\alpha)$ where $\alpha:A\lr A$ and $\mu:A\o A\lr A$ are linear maps,
 with notation $\mu(a\o b)=ab$ such that
 \begin{eqnarray*}
&&\alpha(ab)=\alpha(a)\alpha(b),\ \alpha(1_{A})=1_{A},\\
&&1_{A}a=\alpha(a)=a1_{A},\ \alpha(a)(bc)=(ab)\alpha(c),
\end{eqnarray*}
for any $a,b,c\in A$.
\end{definition}

A linear map $f:(A,\mu_{A},\alpha_{A})\lr (B,\mu_{B},\alpha_{B})$ is called a morphism of Hom-algebra
if $\alpha_{B}\circ f=f\circ\alpha_{A}$, $f(1_{A})=1_{B}$ and $f\circ\mu_{A}=\mu_{B}\circ(f\o f).$

\begin{definition}
  A Hom-Lie algebra
is a triple $(L,[\cdot,\cdot],\alpha)$ consisting of  a linear space  $L$,  a linear map
 $[\cdot,\cdot]: L\o L\rightarrow L$ and  an algebra morphism $\alpha: L\rightarrow L$
 satisfying
\begin{eqnarray*}
&&\alpha[a,a']=[\alpha(a),\alpha(a')],
[a,a']=-[a',a],\\
&&[\a(a),[a',a'']]+[\a(a'),[a'',a]]+[\a(a''),[a,a']]=0,
\end{eqnarray*}
  for any  $a,a',a''\in L$.
\end{definition}
\begin{definition}
Let $(L,[\cdot,\cdot],\alpha)$ be a Hom-Lie algebra. A representation of $L$ is a triple $(V, \rho, \b)$,
where $V$ is a vector space,
$\b\in End(V)$ and $\rho : L \rightarrow End(V)$ is a linear map satisfying
\begin{eqnarray*}
&&\b(\rho(x)v)=\rho(\a(x))\b(v),~\forall x\in L, v\in V,\\
&&\rho([x,y])\circ \b=\rho(\a(x))\circ \rho (y)-\rho(\a(y))\circ \rho(x),
\end{eqnarray*}
 for any $x,y\in L.$
\end{definition}

Let $(M,\mu)$ and $(M',\mu')$ be two left $A$-modules,
then a linear map $f:M\rightarrow N$ is a called left $A$-module map if
\begin{eqnarray*}
 f(a\cdot m)=a\cdot f(m),~f\circ\mu=\mu'\circ f,
\end{eqnarray*}
for any $a\in A$ and $m\in M$.
\medskip

Let $(A, \circ)$ be a vector space with a bilinear operation $\circ: A\o A\rightarrow A$.

(a) Let $L_{\circ}(x)$ and $R_{\circ}(x)$ denote the left and right multiplication operator respectively, that is,
$L_{\circ}(x)y = R_{\circ}(y)x = x\circ y$ for any $x, y \in A$. We also simply denote them by $L(x)$ and $R(x)$
respectively without confusion. In particular, if $(L, [,], \a)$ is a Hom-Lie algebra, we let $ad_{[,]}(x) = ad(x)$
denote the adjoint operator, that is, $ad_{[,]}(x)y = ad(x)y = [x, y]$ for any $x, y \in L$.
\medskip

(b) Let $r=\sum_ia_i\o b_i\in A\o A$. Set
\begin{eqnarray}
r_{12}=\sum_ia_i\o b_i\o 1,~~r_{13}=\sum_ia_i\o1 \o b_i,~~r_{23}=\sum_i1\o a_i\o b_i,
\end{eqnarray}
where $1$ is the unit,  if  $(A, \circ)$ is unital or a symbol playing a similar role to the unit for the
non-unital cases. The operation between two rs is given in an obvious way. For example,
\begin{eqnarray}
&&r_{12}\circ r_{13}=\sum_{i,j}a_i\circ a_j \o b_i\o b_j,~~r_{13}\circ r_{23}=\sum_{i,j}a_i\o  a_j \o b_i\circ b_j,\nonumber\\
&&r_{23}\circ r_{12}=\sum_{i,j} a_j \o a_i\circ b_i\o b_j.
\end{eqnarray}

(c) Let $V$ be a vector space. Let $\tau: V\o V\rightarrow V\o V$ be the exchanging operator defined by
\begin{eqnarray}
\tau(x\o y)=y\o x,
\end{eqnarray}
for any $x,y \in V.$
\medskip

(d) Let $V_1, V_2$ be two vector spaces and  $T: V_1 \rightarrow V_2$ be a linear map. Denoted the dual (linear)
map $T^{\ast}: V_2^{\ast} \rightarrow V_1^{\ast} $  by
\begin{eqnarray}
\langle v_1, T^{\ast}(v_2^{\ast})\rangle=\langle T(v_1), v_2^{\ast}\rangle,~~~\forall x,y \in V.
\end{eqnarray}
for any $x,y \in V.$
\medskip

(e) Let $V$ be a vector space. For any linear map $\rho : A \rightarrow End_k(V)$,
there is a linear map $\rho^{\ast}: A\rightarrow End_k(V^{\ast})$ defined
by
\begin{eqnarray}
 \langle  \rho^{\ast}(x)v^{\ast}, u\rangle= -\langle v^{\ast},  \rho(x)u\rangle,
\end{eqnarray}
 for any $x\in A, u\in V$ and $v^{\ast}\in V^{\ast}.$

\section{Modules and Matched pairs of Hom-Poisson algebras}
\def\theequation{2.\arabic{equation}}
\setcounter{equation} {0} \hskip\parindent

In this section, we give the definition of bimodules of
Hom-Poisson algebras and show that one can obtain the semidirect product Hom-Poisson algebras
 if and only if  is a bimodule of the Hom-Poisson algebra
 and we investigate how to construct a Hom-Poisson algebra which is the direct sum of two Hom-Poisson algebras.

First, we recall from \cite{SB14}   the definition of a matched pair of Hom-Lie algebras as follows.
\begin{definition}
Let $(L_1, [, ]_1, \a_1)$ and $(L_2, [, ]_2, \a_2)$ be two Hom-Lie algebras.
If there are two linear maps $\rho_1 : L_1 \rightarrow End_k(L_2)$
and $\rho_2 : L_2 \rightarrow End_k(L_1)$ which are representations of $L_1$ and $L_2$, satisfying the
following equations:
\begin{eqnarray}
&&\rho_1(\a_1(x))[a,b]_2-[\rho_1(x)a, \a_2(b)]_2- [\a_2(a), \rho_1(x)b]_2\nonumber\\
&&~~~~~~~~~~~~~~~~~~+\rho_1(\rho_2(a)x)\a_2(b)-\rho_1(\rho_2(b)x)\a_2(a)=0,\\
&&\rho_2(\a_2(a))[x,y]_1-[\rho_2(a)x, \a_1(y)]_1- [\a_1(x), \rho_2(a)y]_1\nonumber\\
&&~~~~~~~~~~~~~~~~~~+\rho_2(\rho_1(x)a)\a_1(y)-\rho_2(\rho_1(y)a)\a_1(x)=0,
\end{eqnarray}
for any $x,y \in L_1$ and $a,b\in L_2$. Then $(L_1, L_2, \rho_1, \rho_2, \a_1,\a_2)$ is called a matched pair of Hom-Lie algebras.
In this case, there exists a Hom-Lie algebra structure on the vector space $(L_1\oplus L_2,\a_1+\a_2) $ given by
\begin{eqnarray}
[x+a,y+b]=[x,y]+\rho_2(a)y-\rho_2(b)x+[a, b]+\rho_1(x)b-\rho_1(y)a,\\
(\a_1+\a_2)(x+a)=\a_1(x)+\a_2(a).
\end{eqnarray}
It is denoted by $(L_1\bowtie L_2, \a_1+\a_2)$.
\end{definition}
Let $(A, \circ, \b)$ be a commutative Hom-associative algebra. Recall that a representation (or a module) of $A$
is a vector space $(V, \nu)$ and an $k$-linear map $\mu: A\rightarrow End_k(V)$ such that $\mu(x\circ y)\nu=\mu(\b(x)) \mu(y)$ for any $x,y\in A$.

\begin{theorem}
 Let  $(A_1, \circ_1, \b_1)$ and $(A_2, \circ_2, \b_2)$  be two commutative Hom-associative algebras.
 If there are linear maps
$\mu_1: A_1\rightarrow End_k(A_2)$ and $\mu_2: A_2\rightarrow End_k(A_2)$ which are representation of $A_1$ and $A_1$,
and they satisfy the following equations
\begin{eqnarray}
\b_2(\mu_1(x)a)=\mu_1(\b_1(x))\b_2(a),~~~~\b_1(\mu_2(a)x)=\mu_2(\b_2(a))\b_1(x),\\
\mu_1(\b_1(x))(a\circ_2b)=(\mu_1(x)a)\circ_2\b_2(b)+\mu_1(\mu_2(a)x)\b_2(b),\\
\mu_2(\b_2(a))(x\circ_1y)=(\mu_2(a)x)\circ_1\b_1(y)+\mu_2(\mu_1(x)a)\b_1(y).
\end{eqnarray}
for any $x,y\in A_1$ and $a,b\in A_2$.
Then $(A_1, A_2, \mu_1, \mu_2, \b_1,\b_2)$ is called a matched pair of commutative
Hom-associative algebras. In this case, there exists a commutative Hom-associative algebra structure $\circ$ on the
vector space $(A_1\oplus A_2, \b_1+\b_2)$ given by
\begin{eqnarray}
(x+a)\circ(y+b)=x\circ_1 y +\mu_2(a)y+\mu_2(b)x+a\circ_2b+\mu_1(x)b+\mu_1(y)a,\\
(\b_1+\b_2)(x+a)=\b_1(x)+\b_2(a).
\end{eqnarray}
It is denoted by $(A_1\bowtie A_2, \b_1+\b_2)$.

 Moreover, every commutative Hom-associative algebra
which is the direct sum of the underlying vector spaces of two subalgebras can be obtained from a
matched pair of commutative Hom-associative algebras.
\end{theorem}

{\bf  Proof.} To improve readability, we omit the sybscripts in $\circ_i$.
For any $x,y,z\in A_1$ and  $a,b,c\in A_2$, using (2.5), (2.8) (2.9), we get
\begin{eqnarray*}
&&(\b_1+\b_2)((x+a)\circ(y+b))\\
&=& (\b_1+\b_2)(x\circ y +\mu_2(a)y+\mu_2(b)x+a\circ b+\mu_1(x)b+\mu_1(y)a)\\
&=& \b_1(x\circ y +\mu_2(a)y+\mu_2(b)x)+\b_2(a\circ b+\mu_1(x)b+\mu_1(y)a)\\
&=& \b_1(x\circ y)+\mu_2(\b_2(a))\b_1(y)+\mu_2(\b_2(b))\b_1(x)+\b_2(a\circ b)\\
&&+\mu_1(\b_1(x))\b_2(b)+\mu_1(\b_1(y))\b_2(a)\\
&=& (\b_1+\b_2)(x+a)\circ  (\b_1+\b_2)(y+b).
\end{eqnarray*}
In order to prove that  $(A_1\bowtie A_2, \b_1+\b_2)$ is a Hom-associative algebra, we have to check that
\begin{eqnarray*}
(\b_1+\b_2)(x+a)\circ[(y+b)\circ(z+c)]=[(x+a)\circ(y+b)]\circ (\b_1+\b_2)(z+c).
\end{eqnarray*}
In fact, we have
\begin{eqnarray*}
&&(\b_1+\b_2)(x+a)\circ[(y+b)\circ(z+c)]\\
&=& (\b_1+\b_2)(x+a)[y\circ z +\mu_2(b)z+\mu_2(c)y+b\circ c+\mu_1(y)c+\mu_1(z)b]\\
&=& (\b_1(x)+\b_2(a))[y\circ z +\mu_2(b)z+\mu_2(c)y+b\circ c+\mu_1(y)c+\mu_1(z)b]\\
&=&\b_1(x)\circ(y\circ z) +\mu_2(\b_2(a))(y\circ z+\mu_2(b)z+\mu_2(c)y)+\mu_2(b\circ c+\mu_1(y)c\\
&&+\mu_1(z)b)\b_1(x)+\b_2(a)\circ (b\circ c+\mu_1(y)c+\mu_1(z)b)+\mu_1(\b_1(x))(b\circ c\\
&&+\mu_1(y)c+\mu_1(z)b)+\mu_1(y\circ z +\mu_2(b)z+\mu_2(c)y)\b_2(a)\\
&=& (x\circ y+\mu_2(a)y+\mu_2(b)x)\circ \b_1(z)+\mu_2(a\circ b+\mu_1(x)b+\mu_1(y)a)\b_1(z)\\
&&+\mu_2(\b_2(c))(x\circ y+\mu_2(a)y+\mu_2(b)x)+ (a\circ b+\mu_1(x)b+\mu_1(y)a) \b_2(c)\\
&&+\mu_1(x\circ y+\mu_2(a)y+\mu_2(b)x)\b_2(c)+\mu_1(\b_1(z))(a\circ b+\mu_1(x)b+\mu_1(y)a)\\
&=&(x\circ y+\mu_2(a)y+\mu_2(b)x)\circ (\b_1(z)+\b_2(c))\\
&=&[(x+a)\circ(y+b)]\circ (\b_1+\b_2)(z+c),
\end{eqnarray*}
as desired.
It is easy to check that $(A_1\bowtie A_2, \b_1+\b_2)$ is commutative, and this finishes the proof.
\hfill $\square$

\begin{definition}$^{\cite{MY14}}$
A Hom-Poisson algebra is a quadruple $(P, \mu, [\cdot, \cdot], \a)$ consisting of a
linear space $P$, bilinear maps $\mu: P\times P\rightarrow P$ and $[\cdot, \cdot]: P\times P\rightarrow P$,  and a linear space
homomorphism $\a: P\rightarrow P$ satisfying

(1) $(P,\mu,\a)$ is a commutative Hom-associative algebra,

(2) $(P, [\cdot, \cdot],\a)$ is a Hom-Lie algebra,

(3) for all $x,y,z\in P$,
\begin{eqnarray}
[\a(x), \mu(y,z)]=\mu(\a(y), [x,z])+\mu(\a(z),[x,y]).
\end{eqnarray}
\end{definition}
The condition (2.10) expresses the compatibility between the multiplication and the Poisson
bracket. It can be reformulated equivalently as
\begin{eqnarray*}
[\mu(x,y), \a(z)]=\mu([x,z], \a(y))+\mu(\a(x),[y,z]).
\end{eqnarray*}

A homomorphism between two Hom-Poisson algebras is defined as a linear map between  two
Hom-Poisson algebras preserving the corresponding operations.

\begin{definition}
Let $(P, [\cdot, \cdot], \circ, \a)$ be a Hom-Poisson algebra, $(V, \b)$ be a vector space
and $S, T : P \rightarrow  End_k (V)$ be
two $k$-linear maps. Then $(V, S, T, \b)$ is called a representation (or module) of $P$ if $(V, S, \b)$
is a representation of $(P, [\cdot, \cdot], \a)$ and $(V, T, \b)$ is a representation of $(P, \circ, \a$)
and they are compatible in the
sense that
\begin{eqnarray}
\b(S(x)v)=S(\a(x))\b(v),~~~~\b(T(x)v)=T(\a(x))\b(v),\\
S(x\circ y)\beta=T(\a(y))S(x)+T(\a(x))S(y),\\
T([x, y])\beta=S(\a(x))T(y)-T(\a(y))S(x),
\end{eqnarray}
for any $x,y\in P$.
\end{definition}

In the case of Poisson algebras, we can construct semidirect product when given bimodules.
Analogously, we have

\begin{proposition}
Let $(P, [\cdot, \cdot], \circ, \a)$ be a Hom-Poisson algebra. Then $(V, S, T, \b)$ is a module of a Hom-Poisson algebra $(P, [\cdot, \cdot], \circ, \a)$ if and only if  the direct sum of vector spaces $P\oplus V$  with the operations defined by (we still denote
the operations by $[\cdot, \cdot]$ and $\circ$):
\begin{eqnarray}
\{x_1+v_1, x_2+v_2\}=\{x_1,x_2\}+S(x_1)v_2-S(x_2)v_1, ~~\forall x_1, x_2\in P, v_1,v_2\in V\\
(x_1+v_1)\circ (x_2+v_2)=x_1\circ x_2+T(x_1)v_2+T(x_2)v_1, ~~\forall x_1, x_2\in P, v_1,v_2\in V
\end{eqnarray}
and the twist map $\a+\b:  P\oplus V\rightarrow P\oplus V$ defined by
\begin{eqnarray}
(\a+\b)[x, v]=[\a(x), \b(v)], (\a+\b)(x, v)= (\a(x), \b(v))~~\forall x\in P, v\in V
\end{eqnarray}
  is a Hom-Poisson algebra.  We denote it by $P\ltimes V$.
\end{proposition}
{\bf Proof.}  The sufficient condition holds obviously.
Here we just verify the necessary condition.
For any $x_1,x_2,x_3\in P$ and $v_1,v_2,v_3\in V$, we have
\begin{eqnarray*}
&&(\a+\b)((x_1+v_1)\circ (x_2+v_2))\\
&=&(\a+\b)(x_1\circ x_2+T(x_1)v_2+T(x_2)v_1) \\
&=& \a(x_1\circ x_2)+\b(T(x_1)v_2)+\b(T(x_2)v_1)\\
&=& \a(x_1\circ x_2)+T(\a(x_1))\b(v_2)+T(\a(x_2))\b(v_1)\\
&=& (\a+\b)(x_1+v_1)\circ (\a+\b) (x_2+v_2).
\end{eqnarray*}
Now  we check the following equality
\begin{eqnarray*}
(\a+\b)(x_1+v_1)\circ [(x_2+v_2) \circ (x_3+v_3)]=[(x_1+v_1)\circ (x_2+v_2)] \circ (\a+\b)(x_3+v_3).
\end{eqnarray*}
For this, we calculate
\begin{eqnarray*}
&&(\a+\b)(x_1+v_1)\circ [(x_2+v_2) \circ (x_3+v_3)]\\
&=& (\a+\b)(x_1+v_1)\circ (x_2\circ x_3+T(x_2)v_3+T(x_3)v_2)\\
&=& \a(x_1)\circ(x_2\circ x_3)+T(\a(x_1))(T(x_2)v_3+T(x_3)v_2)+T(x_2\circ x_3)\b(v_1)
\end{eqnarray*}
\begin{eqnarray*}
&=& (x_1 \circ x_2)\circ \a(x_3)+ [T(x_1)T(x_2)]\b(v_3)+[T(x_3)T(x_1)]\b(v_2)+[T(x_3)T(x_2)]\b(v_1)\\
&=& (x_1 \circ x_2)\circ \a(x_3)+ [T(x_1 \circ x_2)]\b(v_3)+[T(x_1)T(x_3)]\b(v_2)+[T(x_2 \circ x_3)]\b(v_1)\\
&=& (x_1\circ x_2+T(x_1)v_2+T(x_2)v_1)\circ (\a+\b)(x_3+v_3)\\
&=& [(x_1+v_1)\circ (x_2+v_2)] \circ (\a+\b)(x_3+v_3).
\end{eqnarray*}
Also one may check directly  that
$
(x_1+v_1)\circ (x_2+v_2)=(x_2+v_2)\circ (x_1+v_1).
$
So $(P\ltimes V, \circ, \a+\b)$ is a  commutative Hom-associative algebra.

Next we will check that $(P\ltimes V,[\cdot, \cdot],\a+\b)$ is a Hom-Lie algebra.
In fact, for any $x_1,x_2,x_3\in P$ and $v_1,v_2,v_3\in V$, we have
\begin{eqnarray*}
[x_1+v_1,x_2+v_2] &=& [x_1,x_2]+S(x_1)x_2-S(x_2)v_2\\
&=&-([x_2,x_1+S(x_2)v_2]-S(x_1)v_1)\\
&=& -[x_2+v_2,x_1+v_1].
\end{eqnarray*}
So the anti-symmetry holds. For the Hom-Jacobi identity, we have
\begin{eqnarray*}
 &&[(\a+\b)(x_1+v_1),[x_2+v_2,x_3+v_3]]+[(\a+\b)(x_2+v_2),[x_3+v_3, x_1+v_1]]\\
 &&+[(\a+\b)(x_3+v_3),[x_1+v_1, x_2+v_2]]\\
 &=&[(\a+\b)(x_1+v_1),[x_2,x_3]+S(x_2)v_3-S(x_3)v_2]+[(\a+\b)(x_2+v_2),[x_3,x_1]\\
 &&+S(x_3)v_1-S(x_1)v_3]+[(\a+\b)(x_3+v_3),[x_1,x_2]+S(x_1)v_2-S(x_2)v_1]\\
 &=&[\a(x_1),[x_2,x_3] ]+S(\a(x_1))(S(x_2)v_3-S(x_3)v_2)-S([x_2,x_3])\b(v_1)+[\a(x_2),[x_3,x_1]]\\
 &&+S(\a(x_2))(S(x_3)v_1-S(x_1)v_3)-S([x_3,x_1])\b(v_2)+[\a(x_3),[x_1,x_2] ]\\
 &&+S(\a(x_3))(S(x_1)v_2-S(x_2)v_1)-S([x_1,x_2])\b(v_3)\\
 &=&0,
\end{eqnarray*}
as desired.
Finally, we check the condition
\begin{eqnarray*}
&&[(\a+\b)(x_1+v_1), (x_2+v_2)\circ (x_3+v_3)]\\
&=&(\a+\b)(x_2+v_2)\circ [x_1+v_1,x_3+v_3]+(\a+\b)(x_3+v_3)\circ[x_1+v_1,x_2+v_2].
\end{eqnarray*}
In fact, on the one hand, we have
\begin{eqnarray*}
&&[(\a+\b)(x_1+v_1), (x_2+v_2)\circ (x_3+v_3)]\\
&=& [(\a+\b)(x_1+v_1), x_2\circ x_3+T(x_2)v_3+T(x_3)v_2]\\
&=& [\a(x_1),x_2\circ x_3]+S(\a(x_1))(T(x_2)v_3+T(x_3)v_2)-S(x_2\circ x_3)\b(v_1)\\
&\stackrel{(2.12)}{=}&[\a(x_1),x_2\circ x_3]+[S(x_1)T(x_2)]\b(v_3)+[S(x_1)T(x_3)]\b(v_2)\\
&&-[T(x_3)S(x_2)]\b(v_1)-[T(x_2)S(x_3)]\b(v_1).
\end{eqnarray*}
On the other hand, we have
\begin{eqnarray*}
&& (\a+\b)(x_2+v_2)\circ [x_1+v_1,x_3+v_3]+(\a+\b)(x_3+v_3)\circ[x_1+v_1,x_2+v_2]\\
 &=& (\a+\b)(x_2+v_2)\circ([x_1, x_3]+S(x_1)v_3-S(x_3)v_1)\\
 &&+(\a+\b)(x_3+v_3)\circ([x_1,x_2]+S(x_1)v_2-S(x_2)v_1)\\
 &=& \a(x_1)\circ [x_1, x_3]+T(\a(x_2))(S(x_1)v_3-S(x_3)v_1)+T([x_1,x_3])\b(v_2)\\
 &&+\a(x_3)\circ [x_1, x_2]+T(\a(x_3))(S(x_1)v_2-S(x_2)v_1)+T([x_1,x_2])\b(v_3)\\
&\stackrel{(2.13)}{=}& \a(x_1)\circ [x_1, x_3]+[T(x_2)S(x_1)]\b(v_3)-[T(x_2)S(x_3)]\b(v_1)+ [S(x_1)T(x_3)] \b(v_2)\\
&&-[T(x_3)S(x_1)] \b(v_2)+\a(x_3)\circ [x_1, x_2]+[T(x_3)S(x_1)]\b(v_2)\\
&& -[T(x_3)S(x_2)]\b(v_1)+[S(x_1)T(x_2)] \b(v_3)-[T(x_2)S(x_1)] \b(v_3)\\
&=& [\a(x_1),x_2\circ x_3]+\a(x_3)\circ [x_1, x_2]+[S(x_1)T(x_2)]\b(v_3)+[S(x_1)T(x_3)]\b(v_2)\\
&&-[T(x_3)S(x_2)]\b(v_1)-[T(x_2)S(x_3)]\b(v_1).
\end{eqnarray*}
Thus $(P\ltimes V, \circ, [\cdot, \cdot], \a+\b)$ is a Hom-Poisson algebra.
The proof is completed.\hfill $\square$

\begin{theorem} Let $(P, [\cdot, \cdot], \circ, \a)$ be a Hom-Poisson algebra and $(V, S, T, \b)$  a module.
Define $S^{\ast}, T^{\ast}: P \rightarrow End_k(V^{\ast})$ by
\begin{eqnarray}
<S^{\ast}(x)u^\ast, v>=-<S(x)v, u^\ast>,~<T^{\ast}(x)u^\ast, v>=-<T(x)v, u^\ast>.
\end{eqnarray}
Assume
$\a^{\ast}: P^{\ast}\rightarrow P^{\ast}, \b^{\ast}: V^{\ast}\rightarrow V^{\ast}$ be the dual maps of $\a$ and $\b$  respectively, that is,
\begin{eqnarray}
\a^{\ast}(x^\ast)(y)=x^{\ast}(\a(y)),~\b^{\ast}(u^\ast)(v)=u^{\ast}(\b(v)).
\end{eqnarray}
 If, in addition,
 \begin{eqnarray}
\b(S(\a(x))v)=S(x)\b(v),~~~~\b(T(x)v)=T(x)\b(v),\\
S(x\circ y)\beta=S(x)T(\a(y))+S(y)T(\a(x)),\\
T([x, y])\beta=T(y)S(\a(x))-S(x)T(\a(y)),
\end{eqnarray}
for any $x,y\in P, x^{\ast}\in P^\ast, u^\ast\in V^\ast, v\in V$.  Then $(V^{\ast},  S^{\ast}, -T^{\ast}, \b^{\ast})$ is a module of $P$.
Moreover,  $(P\ltimes V^{\ast}, \circ, [\cdot, \cdot], \a^{\ast}+\b^{\ast})$ is also a Hom-Poisson algebra.

\end{theorem}

{\bf Proof.} For any  $x,y\in P, x^{\ast}\in P^\ast, u^\ast\in V^\ast, v\in V$, according to (2.17) and (2.18), we have
\begin{eqnarray*}
<\b^{\ast}(S^{\ast}(x)u^{\ast}), v>&=& -<S^{\ast}(x)u^\ast, \b(v)>\\
&=& - <S(x)\b(v), u^\ast>\\
&=&-<\b(S(\a(x))v), u^{\ast}>\\
&=&-<S(\a(x))v, \b^{\ast}(u^{\ast})>\\
&=&<S^{\ast}(\a(x))\b^{\ast}(u^{\ast}), v>.
\end{eqnarray*}
So (2.11) holds for $S^\ast$. Similarly, (2.11) holds for $-T^{\ast}$. According to (2.17)-(2.21), we obtain
\begin{eqnarray*}
&&<S^\ast(x\circ y)\beta^\ast(u^\ast)+T^\ast(\a(y))S^\ast(x)u^\ast+T^\ast(\a(x))S^\ast(y)u^\ast, v>\\
&=& <-\b(S(x\circ y)v)+S(x)T(\a(y))+S(y)T(\a(x)), u^{\ast}>
=0,\\
&&<-T^{\ast}([x, y])\beta^{\ast}(u^{\ast})+S^{\ast}(\a(x))T^{\ast}(y)u^{\ast}-T^{\ast}(\a(y))S^{\ast}(x)u^{\ast},v>\\
&=& <-\b(T([x, y])v)+T(y)S(\a(x))v-S(x)T(\a(y))v,u^{\ast}>
=0.
\end{eqnarray*}
Therefore, (2.12) and (2.13) hold for $(S^\ast, -T^\ast, \a^\ast, \b^\ast)$. It follows that $(S^\ast, -T^\ast, \a^\ast, \b^\ast)$ is a module of $(P, [\cdot, \cdot], \circ, \a)$. The remaining results follow from Proposition 2.5   directly. \hfill $\square$

\begin{example}
Let $(P, [\cdot, \cdot], \circ, \a)$ be a Hom-Poisson algebra. Then $(P, ad_{[\c,\c]}, L_{\circ}, \a)$ and $(P^{\ast},ad^{\ast}_{[\c,\c]}, -L^{\ast}_{\circ}, \a^{\ast})$ are modules of a  Hom-Poisson algebra $(P, [\cdot, \cdot], \circ, \a)$.
\end{example}

In the sequel, we will present a method constructing a Hom-Poisson algebra structure on a direct sum
$P_1\oplus P_2$ of the underlying vector spaces of two Hom-Poisson  algebras $P_1$ and $P_2$ such that $P_1$ and $P_2$
are Hom-Poisson subalgebras.

\begin{theorem}
Let $(P_1, [\cdot, \cdot]_1, \circ_1, \a_1)$ and $(P_2, [\cdot, \cdot]_2, \circ_2, \a_2)$ be two Hom-Poisson algebras.
Assume $\rho_{[\cdot, \cdot]_1}, \mu_{\circ_1}: P_1\rightarrow End_k(P_2)$ and $\rho_{[\cdot, \cdot]_2}, \mu_{\circ_2}: P_2\rightarrow End_k(P_1)$
are four  linear maps such that $(P_1, P_2, \rho_{[\cdot, \cdot]_1}, \rho_{[\cdot, \cdot]_2})$ is a matched pair of Hom-lie algeras and $(P_1, P_2, \mu_{\circ_1}, \mu_{\circ_2})$ is a matched pair of commutative Hom-associative algebras. If in addition, $(P_2, \rho_{[\cdot, \cdot]_1}, \mu_{\circ_1})$ and $(P_1, \rho_{[\cdot, \cdot]_2}, \mu_{\circ_2})$ are representations of the Hom-Poisson algebras $(P_1, [\cdot, \cdot]_1, \circ_1, \a_1)$ and $(P_2, [\cdot, \cdot]_2, \circ_2, \a_2)$ respectively, and $\rho_{[\cdot, \cdot]_1}, \rho_{[\cdot, \cdot]_2},\mu_{\circ_1},  \mu_{\circ_2}$ are compatible in the following sense:
\begin{eqnarray}
\rho_{[\cdot, \cdot]_2}(\a_2(a))(x\circ_1 y)=(\rho_{[\cdot, \cdot]_2}(a)x)\circ_1 \a_1(y)+\a_1(x)\circ_1(\rho_{[\cdot, \cdot]_2}(a)y)\nonumber\\
-\mu_{\circ_2}(\rho_{\{,\}_1}(x)a)\a_1(y)-\mu_{\circ_2}(\rho_{[\cdot, \cdot]_1}(y)a)\a_1(x),~~~\\
\{\a_1(x), \mu_{\circ_2}(a)y\}_1-\rho_{[\cdot, \cdot]_2}(\mu_{\circ_1}(y)a)\a_1(x)=\mu_{\circ_2}(\rho_{[\cdot, \cdot]_1}(x)a)\a_1(y)\nonumber\\
-(\rho_{[\cdot, \cdot]_2}(a)x)\circ_1 \a_1(y)+\mu_{\circ_2}(\a_2(a))([x, y]_1),~~~\\
\rho_{[\cdot, \cdot]_1}(\a_1(x))(a\circ_2 b)=(\rho_{[\cdot, \cdot]_1}(x)a)\circ_2 \a_2(b)+\a_2(a)\circ_2(\rho_{[\cdot, \cdot]_1}(x)b)\nonumber\\
-\mu_{\circ_1}(\rho_{[\cdot, \cdot]_2}(a)x)\a_2(b)-\mu_{\circ_1}(\rho_{[\cdot, \cdot]_2}(b)x)\a_2(a),~~~\\
\{\a_2(a), \mu_{\circ_2}(x)b]_2-\rho_{[\cdot, \cdot]_1}(\mu_{\circ_2}(b)x)\a_2(a)=\mu_{\circ_1}(\rho_{[\cdot, \cdot]_1}(a)x)\a_2(b)\nonumber\\
-(\rho_{[\cdot, \cdot]_1}(x)a)\circ_1 \a_2(b)+\mu_{\circ_1}(\a_1(x))([a, b]_2),
\end{eqnarray}
for any $x,y\in P_1$ and $a,b\in P_2$, then there exists a Hom-Poisson algebra structure ($[\cdot, \cdot], \circ$) on the vector space $P_1\oplus P_2$ given by
\begin{eqnarray}
[x+a,y+b]=([x,y]_1+\rho_{[\cdot, \cdot]_2}(a)y-\rho_{[\cdot, \cdot]_2}(b)x)+([a,b]_2+\rho_{[\cdot, \cdot]_1}(x)b-\rho_{[\cdot, \cdot]_1}(y)a),~~~\\
(x+a)\circ(y+b)=(x\circ_1 y+\mu_2(a)y+\mu_2(b)x)+(a\circ_2 b+\mu_1(x)b+\mu_1(y)a),~~~
\end{eqnarray}
where the twist map $\g:  P_1\oplus P_2\rightarrow P_1\oplus P_2$ is defined by
\begin{eqnarray*}
\g[x, a]=[\a_1(x), \a_2(a)], ~~\forall x\in P_1, a\in p_2.
\end{eqnarray*}
We denote this Hom-Poisson algebra by $P_1\vartriangleright\vartriangleleft P_2$. Moreover $(P_1, P_2, \rho_{[\cdot, \cdot]_1},\mu_{\circ_1}, \rho_{[\cdot, \cdot]_2},$\\$\mu_{\circ_2}, \a_1, \a_2)$ satisfying the above conditions is called a matched pair of Hom-Poisson algebras.
\end{theorem}
{\bf Proof.} To improve readability, we omit the sybscripts in $\circ_i$. For any $x,y\in P_1$ and $a,b\in P_2$,
by Theorem 2.2, $(P_1,P_2,\mu_{\circ_1}, \mu_{\circ_2}, \a_1,\a_2)$ is a matched pair of commutative Hom-associative algebras
if and only if
\begin{eqnarray*}
(a\circ x) \circ \a_1(y)=\a_2(a)\circ (x\circ y),~~(x\circ a) \circ \a_1(y)=\a_1(x)\circ (a\circ y), \\
(x\circ a) \circ \a_2(b)=\a_1(x)\circ (a\circ b),~~(a\circ x) \circ \a_2(b)=\a_2(a)\circ (x\circ b).
\end{eqnarray*}
By Theorem 2.3,   $(P_1,P_2,\rho_{[\c,\c]_1}, \rho_{[\c,\c]_2}, \a_1,\a_2)$ is a matched pair of  Hom-Lie algebras if and only if
\begin{eqnarray*}
&&[\a_2(a),[x,y]]+[\a_1(y),[a,x]]+[\a_1(x), [y,a]]=0,\\
&&[\a_1(x),[a,b]]+[\a_2(b),[x,a]]+[\a_2(a), [b,x]]=0.
\end{eqnarray*}
$(\rho_{[\c,\c]_1},\mu_{\circ_1}, \a_1)$ and $(\rho_{[\c,\c]_2},\mu_{\circ_2}, \a_2)$ satisfy the compatible conditions (2.11-(2.13) and (2.22)-(2.25) if and only if
\begin{eqnarray*}
&&[\a_2(a), x\circ y]=[a,x]\circ \a_1(y)+\a_1(x)\circ [a,y],\\
&&[\a_1(x), a\circ y]=[x,a]\circ \a_1(y)+\a_2(a)\circ [x,y],\\
&&[\a_1(x), a\circ b]=[x,a]\circ \a_2(b)+\a_2(a)\circ [x,b],\\
&&[\a_2(a), x\circ b]=[a,x]\circ \a_2(b)+\a_1(x)\circ [a,b].
\end{eqnarray*}
And this completes the proof.
 \hfill $\square$

\section{Manin triples of Hom-Poisson algebras and Hom-Poisson bialgebras}
\def\theequation{3.\arabic{equation}}
\setcounter{equation} {0} \hskip\parindent

In this section, we introduce the notion of Manin triple of Hom-Poisson algebras which is an analogue of the
notion of Manin triple for Hom-Lie algebras.

\begin{definition}
A Manin triple of Hom-Poisson algebras $(P, P^{+}, P^{-}, \a)$ is a triple of Hom-Poisson algebras $(P,[\c,\c], \circ, \a)$,
$P^{+}$, and $P^{-}$  together with a nondegenerate symmetric bilinear form $B(,)$ on $P$ which is invariant
in the sense that
\begin{eqnarray*}
B([x,y], \a(z))=B(\a(x),[y,z] ),~~~~B(x\circ y, \a(z))=B(\a(x),  y\circ z), ~~\forall x,y,z\in P,
\end{eqnarray*}
satisfying the following conditions:

(1) $P^{+}$ and $P^{-}$ are Hom-Poisson subalgebras of $P$;

(2) $P=P^{+}\oplus P^{-}$ as $k$-vector spaces;

(3) $P^{+}$  and $P^{-}$ are isotropic with respect to $B(,)$.
\end{definition}

A homomorphism between two Manin triples of Hom-Poisson algebras $(P_1, P^{+}_1, P^{-}_1, \a_1)$ and
$(P_2, P^{+}_2, P^{-}_2, \a_2)$ with two nondegenerate symmetric invariant bilinear forms $B_1$ and $B_2$ respectively
is a homomorphism of Hom-Poisson algebras $\varphi: P_1\rightarrow P_2$ such that
\begin{eqnarray*}
\varphi\circ \a_1=\a_2\circ \varphi,~~ \varphi(P^{+}_1)\subset P^{+}_2, ~~~ \varphi(P^{-}_1)\subset P^{-}_2, ~~~B_1(x,y)=B_2(\varphi(x),\varphi(y)).
\end{eqnarray*}

Obviously, a Manin triple of Hom-Poisson algebras is just a triple of Hom-Poisson algebras such that they
are both a Manin triple of Hom-Lie algebras and a commutative Hom-associative version of Manin triple
with the same nondegenerate symmetric bilinear form (and share the same isotropic subalgebras).

Let $(P, [\c,\c], \circ, \a)$ be a Hom-Poisson algebra. If there is a Hom-Poisson algebra structure on the direct sum of the
underlying vector space of $P$ and its dual space $P^{\ast}$ such that $(P\oplus P^{\ast}, P, P^{\ast})$ is a Manin triple of
Hom-Poisson algebras with an invariant symmetric bilinear form on $P\oplus P^{\ast}$ given by
\begin{eqnarray}
B(x+a^{\ast}, y+b^{\ast})=<x,b^{\ast}>+<a^{\ast},y>,\\
(\a+\a^{\ast})(x+a^{\ast})=\a(a)+\a^{\ast}(a^\ast),
\end{eqnarray}
 for any $x,y\in P$ and $a^{\ast},b^{\ast}\in P^{\ast}.$ Then $(P\oplus P^{\ast}, P, P^{\ast}, \a+\a^\ast)$ is called a standard Manin triple of Hom-Poisson algebras. Furthermore, it is straightforward to get the following structure.

\begin{theorem}
Let $(P, [\cdot, \cdot]_1, \circ_1, \a)$ and $(P^{\ast}, [\cdot, \cdot]_2, \circ_2, \a^{\ast})$ be two Hom-Poisson algebras. Then
 $(P\oplus P^{\ast}, P, P^{\ast}, \a+\a^{\ast})$   is a standard Manin triple of Hom-Poisson algebras if and only if $(P, P^{\ast}, ad^{\ast}_{[\cdot, \cdot]_1}, -L^{\ast}_{\circ_1}, ad^{\ast}_{[\cdot, \cdot]_2}, -L^{\ast}_{\circ_2}, \a,\a^{\ast})$ is a matched pair of Hom-Poisson algebras.
\end{theorem}

\begin{definition}

Let $(L, [,], \a)$ be a Hom-Lie algebra and $(L, \delta, \a)$  a Hom-Lie coalgebra.
If $\delta$ is a 1-cocycle of $L$ with coefficient in $L\o L$,
i.e., it satisfies the following equation:
\begin{eqnarray}
\delta([x,y])=(ad(x)\o \a+\a\o ad(x))\delta(y)-(ad(y)\o \a+\a\o ad(y))\delta(x),~~\forall x,y\in L.
\end{eqnarray}
Then we call $(L, [,], \delta, \a)$  a Hom-Lie bialgebra.
In particular, a Hom-Lie bialgebra $(L, [,], \delta, \a)$ is called coboundary if $\delta$ is a 1-coboundary, that is, there exists
an $r\in L\o L$ with $\a^{\otimes 2}(r)=r$ such that
\begin{eqnarray}
\delta(x)=(ad(x)\o \a+\a\o ad(x))r,~~\forall x\in L.
\end{eqnarray}
In this case, $\delta$ automatically satisfies Eq. (3.3). It is usually denoted by $(L, [, ], r, \alpha )$.
\end{definition}
\begin{definition}
Let $(P, \circ,  \a)$ be a  Hom-associative algebra. An infinitesimal Hom-bialgebra structure on $P$ is a linear
map $\Delta: P\rightarrow P\o P$ such that $(P, \Delta)$ is a Hom-coassociative coalgebra and $\Delta$ is a derivation (1-cocycle),
that is,
\begin{eqnarray}
\Delta(a\circ b)=(L_{\circ}(\a(a))\o \a)\Delta(b)+(\a\o R_{\circ}(\a(b)))\Delta(a)     ~~\forall a, b\in P.
\end{eqnarray}
In particular, an infinitesimal Hom-bialgebra $(P, \circ, \Delta, \a)$ is called coboundary if $\Delta$ is a principal derivation
(1-coboundary), that is, there exists an $r\in P\o P$ with $\a^{\otimes 2}(r)=r$ such that
\begin{eqnarray}
\Delta(a)=(L_{\circ}(a)\o \a-\a\o R_{\circ}(a))r,~~\forall a\in P.
\end{eqnarray}
In this case, $\Delta$ automatically satisfies Eq. (3.5). It is usually denoted by $(P, \circ, r)$.
\end{definition}
\begin{definition}
Let $(P, \{, \}, \circ, \a)$ be a Hom-Poisson algebra. Suppose that it is equipped with two
comultiplications $\delta, \Delta: P\rightarrow P\o P$ such that $(P, \delta, \Delta, \a)$ is a Hom-Poisson coalgebra, that is, $(P, \delta)$ is a
Hom-Lie coalgebra and $(P, \Delta)$ is a Hom-coassociative coalgebra which is cocommutative and they satisfy the
following compatible condition:
\begin{eqnarray*}
(\a\o \Delta)\delta(x)=(\delta\o \a)\Delta(x)+(\tau\o id)(\a\o \delta)\Delta(x),~\forall x\in P.
\end{eqnarray*}
If in addition, $(P, \{, \}, \delta, \a)$ is a Hom-Lie bialgebra and $(P, \circ, \Delta, \a)$ is a commutative and cocommutative
infinitesimal Hom-bialgebra and $\delta$ and $\Delta$ are compatible in the following sense
\begin{eqnarray}
\delta(x\circ y)&=&(L_{\circ}(\a(y))\o \a)\delta(x)+(L_{\circ}(\a(x))\o \a)\delta(y)\nonumber\\
&&-(\a\o ad_{\{, \}}(x))\Delta(y)-(\a\o ad_{\{, \}}(y))\Delta(x),\\
\Delta(\{x, y\})&=&(ad_{\{, \}}(\a(x))\o \a +\a\o ad_{\{, \}}(\a(x))) \Delta(y)\nonumber\\
&&+(L_{\circ}(\a(y))\o \a-\a \o L_{\circ}(\a(y))  )\delta(x),
\end{eqnarray}
then $(P, \{, \}, \circ, \delta, \Delta, \a)$ is called a Hom-Poisson bialgebra.
\end{definition}
A homomorphism between two Hom-Poisson coalgebras
is defined as an  linear map between the two Hom-Poisson coalgebras that preserves the corresponding
cooperations. A homomorphism between two Hom-Poisson bialgebras is a homomorphism of both Hom-Poisson
algebras and Hom-Poisson coalgebras.

\begin{theorem} Let $(P, [\cdot, \cdot]_1, \circ_1)$ be a Hom-Poisson algebra equipped with two comultiplications
$\delta, \Delta: P\rightarrow P\o P$. Suppose that $\delta^{\ast}, \Delta^{\ast}: P^{\ast}\o P^{\ast}\subset (P\o P)^{\ast}\rightarrow P^{\ast}$ induce a Hom-Poisson algebra structure on $(P^{\ast}, \a^{\ast})$, where $\delta^{\ast}$ and $\Delta^{\ast}$ correspond to the Hom-Lie bracket and the product of the commutative Hom-associative algebra respectively. Set $[\c, \c]_2=\delta^{\ast}$, $\circ_2=\Delta^{\ast}$. Then the following conditions
are equivalent:

(1) $(P, [\cdot, \cdot]_1, \circ_1, \delta, \Delta, \a)$ is a Hom-Poisson bialgebra.

(2)  $(P, P^{\ast}, ad^{\ast}_{[\cdot, \cdot]_1}, -L^{\ast}_{\circ_1}, ad^{\ast}_{[\cdot, \cdot]_2}, -L^{\ast}_{\circ_2})$ is a matched pair of Hom-Poisson algebras.

(3) $(P\oplus P^{\ast}, P, P^{\ast}, \a+\a^{\ast})$  is a standard Manin triple of Hom-Poisson algebras with the bilinear form defined
by Eq. (3.3) and the isotropic subalgebras are $P$ and $P^{\ast}$.
\end{theorem}
{\bf Proof.} We only prove that the fact (1) holds if and only if the fact (2) holds. In fact, by Theorem 2.8, we only need to check that (2.22)-(2.25) is equivalent to (3.7)-(3.8) in the case that $\rho_{[\c,\c]_{i}}=ad^{\ast}_{[\c,\c]_{i}}, \mu_{\circ_i}=-L_{\circ_i}^{\ast}(i=1,2)$. For any $x,y\in P$ and $a^{\ast}\in P^\ast$, we have
\begin{eqnarray*}
0&=&<-ad^{\ast}_{[\c,\c]_2}(\a^{\ast}(a^{\ast}))(x\circ_1 y)+(ad^{\ast}_{[\c,\c]_2}(a^{\ast})x)\circ_1 \a(y)\\
&&+L_{\circ_2}^{\ast}(ad^{\ast}_{[\c,\c]_1}(x)a^{\ast})\a(y)
+\a(x)\circ_1(ad^{\ast}_{[\c,\c]_2}(a^{\ast})y)+L_{\circ_2}^{\ast}(ad^{\ast}_{[\c,\c]_1}(y)a^{\ast})\a(x), b^{\ast}>\\
&=& -<x\circ_1y,[b^{\ast}, \a^{\ast}(a^{\ast})]_2>-<x,[L_{\circ_1}^{\ast}(\a(y))(b^{\ast}), \a^{\ast2}(a^\ast)]_2>-<\a(y), b^{\ast}\circ_2\\
&& (ad^{\ast}_{[\c,\c]_1}(x))a^\ast>-<y, [L_{\circ_1}^{\ast}(\a(x))(b^{\ast}), \a^{\ast2}(a^\ast)]_2>-<\a(x), (ad^{\ast}_{[\c,\c]_1}(y))a^\ast)\circ_2 b^{\ast}>\\
&=&<-\delta(x\circ y)+(L_{\circ}(\a(y))\o \a)\delta(x)+(L_{\circ}(\a(x))\o \a)\delta(y)\\
&&+(\a\o ad_{[\c,\c]}(x))\Delta(y)+(\a\o ad_{[\c,\c]}(y))\Delta(x),b^{\ast}\o \a^{\ast}(a^{\ast})>.
\end{eqnarray*}
The proof of other cases is silimar.  \hfill $\square$

\section{Coboundary Hom-Poisson bialgebras}
\def\theequation{4.\arabic{equation}}
\setcounter{equation} {0} \hskip\parindent

In this section, we recall some classical conclusions of the coboundary Hom-Lie bialgebras   and give some  important results of the coboundary Hom-Poisson bialgebras.

\begin{definition}$^{\cite{Y09}}$
Let $(L, [,],\a)$ be a Hom-Lie algebra and $r\in L\o L$ with $\a^{\otimes 2}(r)=r$. The linear map $\delta$ defined by Eq. (3.4) makes $(L, \delta)$ into a
Hom-Lie coalgebra if and only if the following conditions are satisfied (for any $x\in L$):
\begin{eqnarray*}
&&(1)(ad(x)\o \a+\a\o ad(x))(r+\tau(r))=0,\\
&&(2)(ad(x)\o \a\o \a+\a\o ad(x)\o \a+\a\o \a\o ad(x))([r_{12}, r_{13}]+[r_{12}, r_{23}]+[r_{13}, r_{23}])=0
\end{eqnarray*}
In particular, the following equation
\begin{eqnarray}
C(r)=[r_{12}, r_{13}]+[r_{12}, r_{23}]+[r_{13}, r_{23}]=0
\end{eqnarray}
is called classical Hom-Yang-Baxter equation(CHYBE).
\end{definition}
\begin{definition}$^{\cite{Y10}}$
Let $(P, \circ, \a)$ be a commutative Hom-associative algebra and $r\in P\o P$ with $\a^{\otimes 2}(r)=r$.
 The linear map $\Delta$ defined by
Eq. (3.6) makes $(P, \Delta)$ into a cocommutative Hom-coassociative coalgebra if and only if the following
conditions are satisfied (for any $x\in P$):

$(1)~(L_{\circ}(x)\o \a-\a\o L_{\circ}(x))(r+\tau(r))=0,$

$(2)~(L_{\circ}(x)\o \a\o \a-\a\o \a\o L_{\circ}(x))(r_{13}\circ r_{12}-r_{12}\circ r_{23}+ r_{23}\circ r_{23})=0.$

\noindent In particular, the following equation
\begin{eqnarray}
A(r)=r_{13}\circ r_{12}-r_{12}\circ r_{23}+ r_{23}\circ r_{23}=0
\end{eqnarray}
is called Hom-associative Yang-Baxter equation (HAYBE).
\end{definition}
Next we will introduce the notion of coboundary Hom-Poisson bialgebra.

\begin{definition}
A Hom-Poisson bialgebra $(P, \{, \}, \circ, \delta, \Delta, \a)$ is called coboundary,  that is, there exists
an $r\in P\o P$ with $\a^{\otimes 2}(r)=r$ such that (for
any $x\in P$)
\begin{eqnarray}
\delta(x)=(ad_{[\cdot, \cdot]}(x)\o \a+\a\o ad_{[\cdot, \cdot]})(x))r,\\
\Delta(x)=(L_\circ(x)\o \a-\a\o L_\circ(x))r,
\end{eqnarray}
\end{definition}

Now  we examine when $(P,\delta, \Delta, \a)$ becomes a Hom-Poisson coalgebra,
where $\delta$ and $\Delta$ are defined by  Eqs. (4.3) and (4.4) for some $r\in P\o P$.
Let $P$ be a vector space equipped with two
comultiplications $\delta,\Delta:P\rightarrow P\o P$. Then $(P,\delta, \Delta)$ becomes a Poisson coalgebra for $\delta$ corresponding
to the Hom-Lie cobracket and $\Delta$ corresponding to the coproduct of the cocommutative Hom-coalgebra if and
only if $(P, \delta)$ is a Hom-Lie coalgebra and $(P, \Delta)$ is a cocommutative Hom-coalgebra and
\begin{eqnarray}
W(x)=(\a\o \Delta)\delta(x)-(\delta\o \a) \Delta(x)-(\tau\o \a)(\a\o \delta)\Delta(x)=0,~~~\forall x\in P.
\end{eqnarray}
Let $(P, [\c,\c], \circ)$ be a Hom-Poisson algebra.  Define two comultiplications $\delta, \Delta: P\rightarrow P\o P$ by Eqs. (4.3) and (4.4) for some $r\in P\o P$, respectively. Then we have
\begin{eqnarray*}
W(x)&=&-(ad_{[\c,\c]}(x)\o \a\o \a)A(r)+(\a \o L_{\circ}(x)\o \a-\a\o \a \o L_{\circ}(x))C(r)\\
&&-\sum_{i}[(ad_{[\c,\c]}(a_i)\o \a)(L_{\circ}(x)\o \a-\a\o L_{\circ}(x))(r+\tau(r))]\o b_i.
\end{eqnarray*}

According to the discussion above, we obtain the following conclusion.

\begin{theorem}
Let $(P, [\cdot, \cdot], \circ, \a)$ be a Hom-Poisson algebra  and $r\in P\o P$ with $\a^{\otimes 2}(r)=r$.
Then the comultiplications $\delta$
and $\Delta$ defined by Eqs. (4.3) and (4.4), make $(P, \delta, \Delta)$ into a Hom-Poisson coalgebra.
Moreover,   $\delta, \Delta, \a)$ is a Hom-Poisson bialgebra if and only if the following conditions are satisfied (for any
$x\in P$):\\
(1) $(ad_{[\cdot, \cdot]}(x)\o \a+\a\o ad_{[\cdot, \cdot]}(x))(r\o \tau(r))=(L_{\circ}(x)\o \a+\a\o L_{\circ}(x))(r+\tau(r))=0$;\\
(2) $(L_{\circ}(x)\o \a\o \a-\a\o \a\o L_{\circ}(x))A(r)=0$;\\
(3)$(ad(x)\o \a\o \a+\a\o ad(x)\o \a+\a\o \a\o ad(x))C(r)=0$;\\
(4)$ W(x) = 0$, where $W(x)$ is given by Eq. (4.5) and $r=\sum a_i\o b_i$.
\end{theorem}

Another important consequence of Theorem 4.4 is the following Hom-Poisson algebra analogue of the
Drinfeld double construction.
\begin{theorem}
Let $(P, [\cdot, \cdot]_1, \circ_1, \delta, \Delta, \a)$ be a Hom-Poisson bialgebra. Then there is a canonical coboundary
Hom-Poisson bialgebra structure on $P+P^{\ast}$.
\end{theorem}
{\bf Proof.} Let $\{e_1,e_2,...,e_n\}$ be a basis of $P$ and $\{e^{\ast}_1,e^{\ast}_2,...,e^{\ast}_n\}$   its dual basis.
Take $r=\sum_i e_{i}\o e^{\ast}_{i}$ with $(\a\o \a^{\ast})r=r$.
Suppose that the Hom-Possion algebra structure $([\c,\c], \circ)$ on $P+P^{\ast}$ is given by $\mathcal{PD}(P)=P\bowtie P^{\ast}$, where $([\c,\c]_2, \circ_2)$ is the Hom-Possion algebra structure on $P^{\ast}$ induced by $(\delta^{\ast}, \Delta^{\ast})$. Then for any $x,y\in P$ and $a^{\ast}, b^{\ast}\in P^\ast$, we have
\begin{eqnarray*}
[x,y]=[x,y]_1,~~x\circ y=x\circ_1 y,~~[a^{\ast}, b^{\ast}]=[a^{\ast}, b^{\ast}]_2, a^{\ast}\circ b^{\ast}=a^{\ast}\circ_2 b^{\ast},\\
ad^{\ast}_{[\c,\c]_1}(x)a^{\ast}=[x,a^{\ast}]=-[a^{\ast}, x]=-ad^{\ast}_{[\c,\c]_2}(a^{\ast})x,~~\\
-L_{\circ_1}^{\ast}(x)a^\ast=x\circ a^{\ast}=a^{\ast}\circ x=--L_{\circ_2}^{\ast}(a^\ast)x.
\end{eqnarray*}
It is straightforward to prove that $r$ satisfies CHYBE and HAYBE and
\begin{eqnarray*}
(ad_{[\cdot, \cdot]}(u)\o \a+\a\o ad_{[\cdot, \cdot]}(x))(r\o \tau(r))=(L_{\circ}(x)\o \a+\a\o L_{\circ}(u))(r+\tau(r))=0,
\end{eqnarray*}
for any $u\in \mathcal{PD}(P)$.
So $r$ satisfies conditions in Theorem 4.4. Thus we have
\begin{eqnarray*}
\delta_{\mathcal{PD}}(u)=(ad_{[\cdot, \cdot]}(u)\o \a^{\ast}+\a\o ad_{[\cdot, \cdot]})(u))r,
 ~~\Delta_{\mathcal{PD}}(u)=(L_\circ(u)\o \a^{\ast}-\a\o L_\circ(u))r,
\end{eqnarray*}
which induce a coboundary Poisson bialgebra structure on $\mathcal{PD}(P)$, as desired. The proof is completed.  \hfill $\square$
\medskip

Let $(P, [\cdot, \cdot]_1, \circ_1, \delta, \Delta, \a)$ be a Hom-Poisson bialgebra. With the Hom-Poisson-bialgebra structure given in
Theorem 4.5, $P+P^{\ast}$ is called the Drinfeld classical double of $P$. As in the proof of Theorem 4.5, we
denoted it by $\mathcal{PD}(P)$.

\section{Post-Hom-Poisson algebra and quasitriangular Hom-Poisson bialgebras}
\def\theequation{5.\arabic{equation}}
\setcounter{equation} {0} \hskip\parindent

In this section,  we introduce notions of $\mathcal{O}$-operator on a Hom-Poisson algebra, post-Hom-Poisson algebra
and quasitriangular Hom-Poisson bialgebra, and show that a quasitriangular Hom-Poisson bialgebra naturally gives a post-Hom-Poisson algebra.

\begin{definition} Let $(L, [\c,\c], \a)$ and $[L',[\c,\c]',\a']$ be two Hom-Lie algebras.
 Suppose that $\rho$ is a Hom-Lie algebra homomorphism from $L$
to $Der_k(L')$, where $Der_k(L')$  consists of all the derivations of $L'$. Then $(L',[\c,\c]', \rho, \a')$ is called a $L$-Hom-Lie algebra.

\end{definition}

\begin{definition} Let $(A, \circ, \a)$ and $(R, \c, \b)$ be two commutative Hom-associative algebras.
Assume $\mu: A\rightarrow End_k(R)$ is a  linear map. Then $(R, \c, \mu, \b)$ is called an $A$-module Hom-algebra if
\begin{eqnarray*}
\mu(x\circ y)=\mu(x)\c \mu(y),~~\mu(\a(x))(u\c v)=(\mu(x)u)\c \b(v),~ \b(\mu(a)v)=\mu(\a(a))\b(v),
\end{eqnarray*}
for any $x,y\in A$ and $u,v\in R$.
\end{definition}
\begin{definition}
Let $(P,[\c,\c], \circ, \a)$ and $(V, [\c,\c]_1, \circ_1, \b)$ be two Hom-Poisson algebras.
Assume  $S, T: P\rightarrow End_k(V)$ are two  linear maps such that

(1) $ (V, [\c,\c]_1, S, \b)$ is a $P$-Hom-Lie algebra, where $P$ is seen as a Hom-Lie algebra with respect to the Hom-Lie bracket
$[\c,\c]$.

(2) $(V, \circ_1, T, \b)$ is a $P$-module Hom-algebra, where $P$ is seen as a commutative Hom-associative algebra
with respect to the commutative associative product $\circ$.

(3) $(V, S, T, \b)$ is a module of $P$.

(4) For any $x\in P$ and $u,v\in V$, the following equations hold:
\begin{eqnarray*}
&&\b(S(x)v)=S(\a(x))\b(v),~~~~\b(T(x)v)=T(\a(x))\b(v),\\
&&S(\a(x))(u\circ_1 v)=(S(x)u)\circ_1 \b(v)+\b(u)\circ_1 (S(x)v),\\
&&[\b(u),T(x)v]_1=-(S(x)u)\circ_1 \b(v)+T(\a(x))[u,v]_1.
\end{eqnarray*}
Then $(V, [\c,\c]_1, \circ_1, S, T, \b)$ is called a $P$-module Hom-Poisson algebra.
\end{definition}

\begin{proposition} With notations as above,  $(V, [\c,\c]_1, \circ_1, S, T, \b)$
is a $P$-module Hom-Poisson algebra if and only if the direct sum of vector spaces $(P\oplus V, \a+\b)$
is turned into a Hom-Poisson algebra by defining the operations as follows:
\begin{eqnarray*}
&&(\a+\b)(x,u)=(\a(x),\b(u)),\\
&&[(x,u), (y,v)]'=([x,y], S(x)v-S(y)u+[u,v]_1),\\
&&(x,u)\circ'(y,v)=(x\circ y, T(x)v+T(y)u+u\circ_1 v),
\end{eqnarray*}
for any $x,y\in P$ and $u,v\in V$.
\end{proposition}

{\bf Proof.} Straightly from Proposition 2.5. \hfill $\square$
\medskip

If $(P,[\c,\c], \circ, \a)$ is a Hom-Poisson algebra,  it is easy to check that  $(P,[\c,\c], \circ,$\\
$ ad_{[\c,\c]}, L_{\circ}, \a)$ is a $P$-module Hom-Poisson algebra.

\begin{definition}
A  Post-Hom-Lie algebra is a vector space $A$ with three bilinear operations $([\c,\c], \diamond, \a)$ satisfying the following equations:
\begin{eqnarray*}
&&[x,y]=-[y,x],\\
&&[[x,y], \a(z)]+[[z,x], \a(y)]+[[y,z], \a(x)]=0,\\
&&\a(z)\diamond (y\diamond x) -\a(y)\diamond (z\diamond x)+ (y\diamond z)\diamond \a(x)- (z\diamond y)\diamond \a(x)+[y,z]\diamond \a(x)=0,\\
&&\a(z)\diamond [x,y]-[z\diamond x,\a(y)]-[\a(x), z\diamond y]=0,
\end{eqnarray*}
 for any $x,y,z\in A$.
\end{definition}

\begin{definition}
A  commutative Hom-dendriform trialgebra is a vector space $A$
equipped with three  bilinear operations $(\c, \succ, \a)$ satisfying the following equations:
\begin{eqnarray*}
&&x \c y=y\c x,\\
&&(x\c y)\c \a(z)=\a(x)\c (y\c z),\\
&&(x\succ y+y\succ x+x\c y)\succ \a(z)=\a(x)\succ (y\succ z),\\
&& (x\succ y)\c \a(z)=\a(x)\succ (y\c z),
\end{eqnarray*}
for any $x,y,z\in A$.
\end{definition}
\begin{definition}
A post-Hom-Poisson algebra is a vector space $A$ equipped with five bilinear
operations $([\c,\c], \diamond, \c, \succ, \a)$  such that $(A, [\c,\c], \diamond, \a)$ is a post-Hom-Lie algebra, $(A, \c, \succ, \a)$
 is a commutative Hom-dendriform trialgebra, and satisfy the following compatible conditions:
 \begin{eqnarray}
&&[\a(x), y\c z]=[x,y]\c \a(z)+\a(y)\c [y,z],\\
&&[\a(x), z\succ y]=\a(z)\succ [x,y]-\a(y)\c (z\diamond x),\\
&&\a(x)\diamond(y \c z)=(x\diamond y)\c \a(z)+\a(y)\c (x\diamond z),\\
&&(y\succ z+z\succ y+y\c z)\diamond \a(x)=\a(z)\succ (y\diamond x)+\a(y)\succ (z\diamond x),\\
&& \a(x)\diamond (z\succ y)=\a(z)\succ (x\diamond y)+(x\diamond z-z\diamond x+[x,z] )\succ \a(y),
 \end{eqnarray}
 for any $x,y,z\in A$.
\end{definition}

Let $(A, [\c,\c], \diamond, \c, \succ, \a)$ be a post-Hom-Poisson algebra, it is obvious that $(A, [\c,\c], \c, \a)$ is a Hom-Poisson algebra.
 On the other hand, we have the following conclusion.

\begin{theorem}
Let $(A, [\c,\c], \diamond, \c, \succ, \a)$ be a post-Hom-Poisson algebra. Define two new bilinear operations
\begin{eqnarray*}
\{x,y\}=x\diamond y-y\diamond x+[x,y],~~~x\circ y=x\succ y+y\succ x+x\c y,~~\forall x,y\in A.
\end{eqnarray*}
Then $(A, \{\c,\c\}, \circ, \a)$ becomes a Hom-Poisson algebra. It is called the associated Hom-Poisson algebra of
$(A, [\c,\c], \diamond, \c, \succ, \a)$ and is denote by $(P(A),\{\c,\c\}, \circ, \a )$.
\end{theorem}
{\bf Proof.} For any $x,y,z\in A$, we will check that the equality
$
\a(x)\circ(y\circ z) =(x\circ y)\circ \a(z).
$
In fact, we have
\begin{eqnarray*}
&&\a(x)\circ(y\circ z)\\
&=& \a(x)\circ(y\succ z+z\succ y+y\c z)\\
&=& \a(x)\succ (y\succ z+z\succ y+y\c z)+(y\succ z+z\succ y\\
&&+y\c z)\succ \a(x)+ \a(x)\c (y\succ z+z\succ y+y\c z)
\end{eqnarray*}
\begin{eqnarray*}
&=& (x\succ y+y\succ x+x\c y)\succ \a(z)+\a(z)\succ (x\succ y\\
&&+y\succ x+x\c y)+ (x\succ y+y\succ x+x\c y)\c \a(z)\\
&=& (x\succ y+y\succ x+x\c y)\circ \a(z)\\
&=& (x\circ y)\circ \a(z).
\end{eqnarray*}
Also one may check directly that
$
x\circ y=y\circ x.
$
So $(A, \circ, \a )$ is a commutative Hom-associative algebra.

Next, we check that $(A, \{\c,\c\}, \a )$ is a  Hom-Lie algebra. For this, we take $x,y,z\in A$ and calculate
\begin{eqnarray*}
\{x,y\}=x\diamond y-y\diamond x+[x,y]
= y\diamond x-x\diamond y-[x,y]
= -\{y,x\}.
\end{eqnarray*}
So the anti-symmetry holds. For Hom-Jacobi identity, we have
\begin{eqnarray*}
&&\{\a(x), \{y,z\}\}+\{\a(y),\{z, x\}+\{\a(z), \{x,y\}\}\}\\
&=&  \{\a(x), y\diamond z-z\diamond y+[y,z]\}+\{\a(y),z\diamond x-x\diamond z\\
&&+[z,x]\}+\{\a(z), x\diamond y-y\diamond x+[x,y]\}\\
&=& \a(x)\diamond (y\diamond z-z\diamond y+[y,z])-(y\diamond z-z\diamond y+[y,z])\diamond \a(x)+[\a(x), [y,z]]\\
&& \a(y)\diamond(z\diamond x-x\diamond z+[z,x])-(z\diamond x-x\diamond z+[z,x])\diamond \a(y)+[\a(y),[z,x]]\\
&& \a(z)\diamond (x\diamond y-y\diamond x+[x,y])-(x\diamond y-y\diamond x+[x,y])\diamond\a(z)+[\a(z),[x,y]]\\
&=&0.
\end{eqnarray*}

Finally, we check the condition
$
\{\a(x), y\circ z\}=\a(y)\circ \{x,z\}+\a(z)\circ \{x,y\}.
$
In fact, we have
\begin{eqnarray*}
&&\{\a(x), y\circ z\}\\
&=& \a(x)\diamond (y \circ z)-(y \circ z)\diamond \a(x)+[\a(x), y\circ z]\\
&=& \a(x)\diamond (y\succ z+z\succ y+ y\c z)-(y\succ z+z\succ y+ y\c z)\diamond \a(x)\\
&&+[\a(x), y\succ z+z\succ y+ y\c z]\\
&=& \a(x)\diamond (y\succ z)+ \a(x)\diamond (z\succ y)+  \a(x)\diamond (y\c z)-\a(z)\succ(y\diamond x)\\
&&-\a(y)\succ(z\diamond x)+[\a(x), y\succ z]+[\a(x), z\succ y]+[\a(x),  y\c z]\\
&=& \a(y)\succ (x\diamond z)+(x\diamond y-y\diamond x +[x,y])\succ \a(z)\\
&&  \a(z)\succ (x\diamond  y)+(x\diamond z-z\diamond x+[x,z])\succ \a(y)\\
&& (x\diamond y)\c \a(z)+\a(y)\c(x\diamond z)-\a(z)\succ(y\diamond x)-\a(y)\succ(z\diamond x)\\
&&\a(y)\succ [x,z]-\a(z)\c (y \diamond x)+\a(z)\succ [x,y]-\a(y)\c (z \diamond x)+[\a(x),  y\c z]
\end{eqnarray*}
\begin{eqnarray*}
&=&\a(y)\succ(x\diamond z-z\diamond x+[x,z])+(x\diamond z-z\diamond x\\
&& +[x,z])\succ \a(y)+\a(y)\c(x\diamond z-z\diamond x+[x,z])\\
 &&\a(z)\succ(x\diamond y-y\diamond x+[x,y])+(x\diamond y-y\diamond x\\
&& +[x,y])\succ \a(z)+\a(z)\c(x\diamond y-y\diamond x+[x,y])\\
&=& \a(y)\succ \{x,z\}+\{x,z\}\succ \a(y)+\a(y)\c \{x,z\}+\a(z)\succ \{x,y\}\\
&& +\{x,y\}\succ \a(z)+\a(z)\c \{x,y\}\\
&=&\a(y)\circ \{x,z\}+\a(z)\circ \{x,y\},
\end{eqnarray*}
as desired.
 Thus $(A, \{\c,\c\}, \circ, \a)$ is a Hom-Poisson algebra. The proof is finished. \hfill $\square$

\begin{definition}
Let $(P, [\c,\c], \circ, \a)$ be a Hom-Poisson algebra and $(V, [\c,\c]_1, \circ_1, S, T, \b)$ a $P$-module Hom-Poisson algebra.
 A linear map $R: V\rightarrow P$ is called an $O$-operator of weight $\lambda\in k$ associated to
$(V, [\c,\c]_1, \circ_1, S, T, \b)$ if, for any $u,v\in V$,
\begin{eqnarray}
&&\a \circ R=R\circ \b,\\
&&[R(u), R(v)]=R(S(R(u))v-S(R(v))u+\lambda[u,v]_1),\\
&&R(u)\circ R(v)=R(T(R(u))v+T(R(v))u+\lambda u\circ_1 v).
\end{eqnarray}
When $(V, [\c,\c]_1, \circ_1, S, T, \b)=(P, [\c,\c], \circ, ad_{[\c,\c]}, L_{\circ}, \a)$, Eqs(5.7)-(5.8) become
\begin{eqnarray}
&&[R(u), R(v)]=R([R(u), v]-[u, R(v)]+\lambda[u,v]),\\
&&R(u)\circ R(v)=R(R(u)\circ v+u \circ R(v)+\lambda u\circ v).
\end{eqnarray}
respectively. Eqs. (5.9) and (5.10) imply that $R: P\rightarrow P$ is a Rota-Baxter operator of weight $\lambda\in k$
on the Hom-Lie algebra $(P, [\c,\c], \a)$ and on the (commutative) Hom-associative algebra $(P, \circ, \a)$ respectively.
\end{definition}

\begin{theorem}
Let $(P, [\c,\c], \circ, \a)$ be a Hom-Poisson algebra and $(V, [\c,\c]_1, \circ_1, S, T, \b)$  a $P$-module Hom-Poisson algebra.
The linear map $R: V\rightarrow P$  is  an $\mathcal{O}$-operator of weight $\lambda\in k$ associated to
$(V, [\c,\c]_1, \circ_1, S, T, \b)$. Define four new bilinear operations $\{\c,\c\}, \diamond, \c, \succ: V\o V\rightarrow V$ as follows:
\begin{eqnarray*}
\{u,v\}=\lambda[u,v]_1,~~u\diamond v=S(R(u))v,~~u\c v=\lambda u\circ_1 v,~~u\succ v=T(R(u))v,
\end{eqnarray*}
for any $u,v\in V$. Then $(V, \{\c,\c\}, \diamond, \c, \succ, \b)$ becomes a post-Hom-Poisson algebra and $R$ is a homomorphism of Hom-Poisson algebras
from the associated Hom-Poisson algebra $P(V)$ of $(V, \{\c,\c\}, \diamond, \c, \succ, \b)$ to $(P, [\c,\c], \circ, \a)$.
\end{theorem}
{\bf Proof.} First we check that $(V, \{\c,\c\}, \diamond, \b)$ is a Post-Hom-Lie algebra, for any $u,v,w\in V$, it is easy to obtain that
\begin{eqnarray*}
&&\{u,v\}=-\{v,u\},\\
&&\{\{u,v\}, \b(w)\}+\{\{w,u\}, \b(v)\}+\{\{v,w\}, \b(u)\}=0.
\end{eqnarray*}
So it is sufficient to verify the following conditions:
\begin{eqnarray*}
&&\b(w)\diamond (v\diamond u) -\b(v)\diamond (w\diamond u)+ (v\diamond w)\diamond \b(u)- (w\diamond v)\diamond \b(u)+\{v,w\}\diamond \b(u)=0,\\
&&\b(w)\diamond \{u,v\}-\{w\diamond u,\b(v)\}-\{\b(u), w\diamond v\}=0.
\end{eqnarray*}
In fact, we have
\begin{eqnarray*}
&&\b(w)\diamond (v\diamond u) -\b(v)\diamond (w\diamond u)+ (v\diamond w)\diamond \b(u)- (w\diamond v)\diamond \b(u)+\{v,w\}\diamond \b(u)\\
&=&S(R(\b(w))) S(R(v))u-S(R(\b(v))) S(R(w))u+S(R(S(R(v))w))\b(u)\\
&&-S(R(S(R(w))v))\b(u)+S(R(\lambda[v,w]_1))\b(u)\\
&=& S(R(w))\circ S(R(v))\b(u)-S(R(v))\circ S(R(w))\b(u)+S([R(v),R(w)])\b(u)\\
&=& 0
\end{eqnarray*}
and
\begin{eqnarray*}
 &&\b(w)\diamond \{u,v\}-\{w\diamond u,\b(v)\}-\{\b(u), w\diamond v\}\\
 &=& \lambda(\b(w)\diamond [u,v]_1-[w\diamond u,\b(v)]_1-[\b(u), w\diamond v]_1)\\
 &=& \lambda(S(R(\b(w)))[u,v]_1-[S(R(w))u, \b(v)]_1-[\b(u), S(R(w))v]_1)\\
 &=& 0,
\end{eqnarray*}
as required.

Next  we will check that $(V, \c, \succ, \b)$ is  a commutative Hom-dendriform trialgebra, for any $u,v,w\in V$, it is easy to obtain that
$u \c v=v\c u$ and
$(u\c v)\c \b(w)=\b(u)\c (v\c w).
$
So we need to check the following conditions:
\begin{eqnarray*}
&&(u\succ v+v\succ u+u\c v)\succ \b(w)=\b(u)\succ (v\succ w),\\
&& (u\succ v)\c \b(w)=\b(u)\succ (v\c w).
\end{eqnarray*}
For the above equality, we calculate
\begin{eqnarray*}
&&(u\succ v+v\succ u+u\c v)\succ \b(w)\\
&=& (T(R(u))v+T(R(v))u+\lambda u\circ_1 v)\succ \b(w)\\
&=& T(R(T(R(u))v+T(R(v))u+\lambda u\circ_1 v))\b(w)\\
&=&( T(R(u))\circ T(R(v)))\b(w)\\
&=& ( T(R(\b(u)))(T(R(v))w)\\
&=& \b(u)\succ (v\succ w).
\end{eqnarray*}
For the below equality, we have
\begin{eqnarray*}
(u\succ v)\c \b(w)
=\lambda T(R(u)v)\circ_1 \b(w)
= T(R(\b(u))( v\c w)
= \b(u)\c (v\c w).
\end{eqnarray*}
The Proof of other cases is similar. Then $(V, \{\c,\c\}, \diamond, \c, \succ, \b)$ becomes a post-Hom-Poisson algebra and
$R$ is a homomorphism of Hom-Poisson algebras
from the associated Hom-Poisson algebra $P(V)$ of $(V, \{\c,\c\}, \diamond, \c, \succ, \b)$ to $(P, [\c,\c], \circ, \a)$.    \hfill $\square$

\begin{definition}
A coboundary Hom-Poisson bialgebra $(P, [\c,\c], \circ, r, \a)$ is called quasitriangular if $r$ is a solution of HPYBE,
that is, it satisfies both HAYBE and CHYBE.
\end{definition}

 Obviously, a coboundary Hom-Poisson bialgebra $(P, [\c,\c], \circ, r, \a)$ is quasitriangular if and
only if $(P, [\c,\c], r, \a)$ is a quasitriangular Hom-Lie bialgebra and $(P,  \circ, r, \a)$ is a quasitriangular infinitesimal Hom-bialgebra.

Let $(V,\b)$ be a finite-dimensional vector space and $r\in V\o V$ such that $(\b\o \b)(r)=r$.
Then $r$ can be identified as a linear map from $V^\ast$ to $V$ which we still denote by $r$ through
\begin{eqnarray*}
 <r(a^{\ast}), b^{\ast}>=<a^{\ast}\o b^{\ast}, r>,~~~\b\circ r=r\circ \b^{\ast},~~\forall a^{\ast},b^{\ast}\in V^\ast.
\end{eqnarray*}
Define a linear map $r':V^{\ast}\rightarrow V$ by
\begin{eqnarray*}
 <a^{\ast},r'(b^{\ast})>=<r, a^{\ast}\o b^{\ast},>,~~~\b\circ r'=r'\circ \b^{\ast}.
\end{eqnarray*}
We call
\begin{eqnarray}
\varphi=\frac{r-r'}{2},~~~~~\psi=\frac{r+r'}{2},
\end{eqnarray}
the skew-symmetric part and the symmetric part of $r$ respectively. Therefore, we have the following theorem.
\begin{theorem}
Let $(P, [\c,\c], \circ, \a, r)$ be a quasitriangular Hom-Poisson bialgebra with the linear map $\psi$ defined by (5.11).
 Define two new bilinear operations $\{\c,\c\}:P\o P\rightarrow P$ as follows:
\begin{eqnarray}
\{a^\ast,b^{\ast}\}=-2ad^{\ast}_{[\c,\c]}(\psi(a^{\ast}))b^{\ast},~~a^{\ast}\c b^{\ast}=2L_{\circ}^{\ast}(\psi(a^{\ast}))b^{\ast}, ~~\forall a^{\ast},b^{\ast}\in P^\ast.
\end{eqnarray}
Then $(P^{\ast},\{\c,\c\}, \c, ad^{\ast}_{[\c,\c]}, -L_{\circ}^{\ast}, \a^{\ast})$ becomes a $P$-module Hom-Poisson algebra of $(P, [\c,\c], \circ, \a)$. Moreover, $r$ regarded as a linear map from $P^\ast$ to $P$ is  an $\mathcal{O}$-operator of weight 1 associated to $(P^{\ast},\{\c,\c\}, \c, ad^{\ast}_{\{\c,\c\}},-L_{\circ}^{\ast}, \a^{\ast})$, that is,
\begin{eqnarray*}
&&\a \circ R=R\circ \a^{\ast},\\
&&[r(a^{\ast}), r(b^{\ast})]=r(ad^{\ast}_{[\c,\c]}(r(a^{\ast}))b^{\ast}-ad^{\ast}_{[\c,\c]}(r(b^{\ast}))a^{\ast}+\{a^{\ast},b^{\ast}\}),\\
&&r(a^{\ast})\circ r(b^{\ast})=r(-L_{\circ}^{\ast}(R(a^{\ast}))b^{\ast}+-L_{\circ}^{\ast}(r(b^{\ast}))a^{\ast}+ a^{\ast}\c b^{\ast}).
\end{eqnarray*}
\end {theorem}

The following result establishes a close relation between a post-Hom-Poisson algebra and a quasitriangular
Hom-Poisson bialgebra.

\begin{theorem}
With the notations above, we define four new bilinear operations $\{\c,\c\}, \diamond, \c, \succ: P^{\ast}\o P^{\ast}\rightarrow P^{\ast}$ as follows:
\begin{eqnarray*}
&&\{a^\ast,b^{\ast}\}=-2ad^{\ast}_{[\c,\c]}(\psi(a^{\ast}))b^{\ast},~~~~ a^{\ast}\diamond b^{\ast}=ad^{\ast}_{[\c,\c]}(r(a^{\ast}))b^{\ast},\\
&&a^{\ast}\c b^{\ast}=2L_{\circ}^{\ast}(\psi(a^{\ast}))b^{\ast}, ~~~a^{\ast}\succ  b^{\ast}=-L_{\circ}^{\ast}(r(a^{\ast}))b^{\ast}.
\end{eqnarray*}
for any $a^{\ast},b^{\ast}\in P^\ast$.
Then $(P^{\ast},\{\c,\c\}, \diamond, \c, \succ, \a^{\ast} )$ becomes a post-Hom-Poisson algebra and $r$ is a homomorphism of Hom-Poisson
algebras from the associated Hom-Poisson algebra $P(P^{\ast})$ of $(P^{\ast},\{\c,\c\}, \diamond, \c, \succ, \a^{\ast} )$ to $(P,[\c,\c], \circ, \a)$.
\end{theorem}

{\bf Proof.} It follows from Theorem 5.10 and 5.12 directly.
\section*{ACKNOWLEDGEMENTS}

The work of S. J. Guo is  supported by  the NSF of China (No. 11761017) and
   the Fund of Science and Technology Department of Guizhou Province (No. [2019]1021).
\medskip

 The work of X. H. Zhang is  supported by   the Project Funded by China Postdoctoral Science Foundation (No. 2018M630768)and
the NSF of Shandong Province (No. ZR2016AQ03).
\medskip

   The work of S. X. Wang is  supported by  the outstanding top-notch talent cultivation project of Anhui Province (No. gxfx2017123)
 and  the NSF of Anhui Provincial  (1808085MA14).


\begin{thebibliography}{aa}


\bibitem{V94}   Vaisman, I. (1994). Lectures on the Geometry of Poisson Manifolds, Progress in Mathematics Vol. 118 (Birkh¡§auser Verlag, Basel).

\bibitem{B10} Bai, C. M. (2010). Double constructions of Frobenius algebras, Connes cocycles and their duality. J. Noncommu. Geom. 4: 475-530.


\bibitem{BGN10} Bai, C. M.,     Guo, L.,   Ni, X. (2010).  Nonabelian generalized Lax pairs, the classical Yang-Baxter equation and Post Lie algebras.
Commun. Math. Phys. 297: 553-596.

\bibitem{BGN12}  Bai, C. M.,  Guo, L.,   Ni, X. (2012). $\mathcal{O}$-operators on associative algebras and associative Yang-Baxter equations.  Pac. J. Math. 256: 257-289.

\bibitem{BM14} Benayadi, S., Makhlouf, A. (2014). Hom-Lie algebras with symmetric invariant nondegenerate bilinear forms.
J. Geom. Phys. 76:38-60.



\bibitem{BEM12} Bordemann, M., Elchinger, O., Makhlouf, A. (2012). Twisting Poisson algebras, CoPoisson algebras and Quantization. Trav. Math. 20: 83-119.

\bibitem{CWZ10} Chen, Y., Wang, Z. W., Zhang, L. Y. (2010). The construction of Hom-Lie bialgebra.   J. Lie Theory   20:767-783.

\bibitem{HLS06}  Hartwig, J., Larsson, D., Silvestrov, S. (2006). Deformations of Lie algebras using $\sigma$-derivations. J. Algebra 295:
314-361.



\bibitem{MS08}  Makhlouf, A., Silvestrov, S. (2008). Hom-algebra structures. J. Gen. Lie Theory Appl. 2(2):51-64.

\bibitem{MY14}  Makhlouf, A.,  Yau, D.  (2014).   Rota-Baxter Hom-Lie-admissible algebras.  Comm. Algebra 42(3): 1231-1257.

\bibitem{NB13} Ni, X.,  Bai, C. M. (2013). Poisson bialgebras.  J. Math. Phys.    54: 023515.

\bibitem{SC13} Sheng, Y. H., Chen, D. (2013). Hom-Lie 2-algebras. J. Algebra 376:174-195.

\bibitem{SB14} Sheng, Y. H., Bai, C. M. (2014). A new approach to Hom-Lie bialgebras. J. Algebra 399:232-250.

\bibitem{Y09} Yau, D. (2009). The Hom-Yang-Baxter equation, Hom-Lie algebras, and quasi-triangular bialgebras. J. Phys. A 42:
165-202.
\bibitem{Y10}  Yau, D. (2010).  Infinitesimal Hom-bialgebras and Hom-Lie bialgebras.  arXiv:1001.5000.


\bibitem{Y11} Yau, D. (2011). The Hom-Yang-Baxter equation and Hom-Lie algebras. J. Math. Phys. 52:053502.
\bibitem{Y15} Yau, D. (2015). The classical Hom-Yang-Baxter equation and Hom-Lie bialgebras. Int. Electron. J. Algebra 17:11-45.



\end{thebibliography}
\end{document}